\RequirePackage{fix-cm}

\documentclass[smallcondensed,nospthms]{svjour3}
\smartqed  
\usepackage{graphicx}
\usepackage{mathptmx}      
\usepackage{amssymb}
\usepackage{amsmath}
\usepackage[all,cmtip]{xy}
\usepackage{amsthm}
\usepackage{color}
\usepackage{enumitem}
\usepackage{stmaryrd}
\usepackage[misc]{ifsym}
\usepackage[colorlinks=true]{hyperref}
\hypersetup{allcolors=[rgb]{0.1,0.1,0.4}}
\input diagxy

\theoremstyle{plain}
  \newtheorem{thm}{Theorem}[section]
  \newtheorem{lem}[thm]{Lemma}
  \newtheorem{prop}[thm]{Proposition}
  \newtheorem{cor}[thm]{Corollary}
\theoremstyle{definition}
  \newtheorem{defn}[thm]{Definition}
  \newtheorem{exmp}[thm]{Example}
  \newtheorem{rem}[thm]{Remark}

\DeclareMathAlphabet{\mathcal}{OMS}{cmsy}{m}{n}

\def\oto{{\bfig\morphism<180,0>[\mkern-4mu`\mkern-4mu;]\place(86,0)[\circ]\efig}}

\newcommand{\da}{\downarrow}
\newcommand{\ua}{\uparrow}
\newcommand{\ra}{\rightarrow}

\newcommand{\lda}{\swarrow}
\newcommand{\rda}{\searrow}

\newcommand{\rat}{\rightarrowtail}
\newcommand{\bv}{\bigvee}
\newcommand{\bw}{\bigwedge}

\newcommand{\dv}{\dashv}

\renewcommand{\phi}{\varphi}
\newcommand{\al}{\alpha}
\newcommand{\be}{\beta}

\newcommand{\Lam}{\Lambda}
\newcommand{\lam}{\lambda}

\newcommand{\CQ}{\mathcal{Q}}

\newcommand{\sC}{{\sf C}}

\newcommand{\sK}{{\sf K}}

\newcommand{\sM}{{\sf M}}

\newcommand{\sV}{{\sf V}}

\newcommand{\sy}{\mathfrak{y}}

\newcommand{\Fix}{{\sf Fix}}

\newcommand{\Cat}{{\bf Cat}}

\newcommand{\VCat}{\sV\text{-}\Cat}

\newcommand{\dphi}{\phi^{\da}}
\newcommand{\uphi}{\phi_{\ua}}
\newcommand{\dphiop}{(\phi^{\op})^{\da}}
\newcommand{\uphiop}{(\phi^{\op})_{\ua}}

\newcommand{\dtau}{\tau^{\da}}
\newcommand{\utau}{\tau_{\ua}}

\newcommand{\syd}{\sy^{\dag}}

\newcommand{\iotad}{\iota^{\dag}}
\newcommand{\co}{{\rm co}}
\newcommand{\op}{{\rm op}}

\newcommand{\Mphi}{\sM\phi}

\newcommand{\Kphi}{\sK\phi}
\newcommand{\Kdphi}{\sK^{\dag}\phi}
\newcommand{\Kpsi}{\sK\psi}
\newcommand{\Kdzeta}{\sK^{\dag}\zeta}

\newcommand{\with}{\mathrel{\circledast}}
\newcommand{\ts}{\mathrel{\otimes}}

\renewcommand{\leq}{\leqslant}
\renewcommand{\geq}{\geqslant}

\newcommand{\LSC}{{\sf LSC}(X,[0,\infty])}

\numberwithin{equation}{section}

\allowdisplaybreaks

\begin{document}

\title{Isbell adjunctions and Kan adjunctions via quantale-enriched two-variable adjunctions\thanks{This is a post-peer-review, pre-copyedit version of an article published in \emph{Applied Categorical Structures}. The final authenticated version is available online at: \url{https://doi.org/10.1007/s10485-021-09654-w}.}
}


\author{Lili Shen         \and
        Xiaoye Tang 
}


\institute{L. Shen \Letter \at
              School of Mathematics, Sichuan University, Chengdu 610064, China\\
              \email{shenlili@scu.edu.cn}           
           \and
           X. Tang \at
              School of Mathematics, Sichuan University, Chengdu 610064, China\\
              \email{tang.xiaoye@qq.com}
}

\date{}

\maketitle

\begin{abstract}
It is shown that every two-variable adjunction in categories enriched in a commutative quantale serves as a base for constructing Isbell adjunctions between functor categories, and Kan adjunctions are precisely Isbell adjunctions constructed from suitable associated two-variable adjunctions. Representation theorems are established for fixed points of these adjunctions.
\end{abstract}

\section{Introduction}

A distributor \cite{Benabou1973} $\phi\colon X\oto Y$ between categories enriched in a commutative quantale \cite{Lawvere1973,Kelly1982,Rosenthal1990,Hofmann2014}
$$\sV=(\sV,\otimes,k)$$
may be described as a $\sV$-bifunctor $\phi\colon X^{\op}\ts Y\to\sV$, and it induces three pairs of adjoint $\sV$-functors between the (co)presheaf $\sV$-categories of $X$ and $Y$:
\begin{enumerate}[label=(\arabic*)]
\item the \emph{Isbell adjunction} \cite{Shen2013a} $\uphi\dv\dphi:(\sV^Y)^{\op}\to\sV^{X^{\op}}$,
\item the \emph{Kan adjunctions} \cite{Shen2013a,Shen2014} $\phi^*\dv\phi_*:\sV^{X^{\op}}\to\sV^{Y^{\op}}$ and $\phi_{\dag}\dv\phi^{\dag}:(\sV^X)^{\op}\to(\sV^Y)^{\op}$.
\end{enumerate}

In this paper we demonstrate that every \emph{two-variable adjunction} \cite{Gray1980,Shulman2006} between $\sV$-categories gives rise to (generalized) Isbell adjunctions and Kan adjunctions, and in fact provides a unified framework for these notions. Explicitly, for every two-variable adjunction $(X,Y,Z,\with,\lda,\rda)$ between $\sV$-categories (Definition \ref{two-var-adj-def}), a $\sV$-bifunctor
$$\phi\colon A^{\op}\ts B\to Z$$
induces an Isbell adjunction
$$\uphi\dv\dphi\colon (X^B)^{\op}\to Y^{A^{\op}}$$
between the $\sV$-categories of $\sV$-functors (Proposition \ref{uphi-dv-dphi}). As every two-variable adjunction is associated with several others (Lemma \ref{two-var-adj-dual}), Isbell adjunctions constructed upon suitable two-variable adjunctions are precisely Kan adjunctions
$$\psi^*\dv\psi_*\colon Z^{A^{\op}}\to X^{B^{\op}}\quad\text{and}\quad\zeta_{\dag}\dv\zeta^{\dag}\colon (Y^A)^{\op}\to(Z^B)^{\op}$$
induced by $\sV$-bifunctors
$$\psi\colon A^{\op}\ts B\to Y\quad\text{and}\quad\zeta\colon A^{\op}\ts B\to X,$$
respectively (Proposition \ref{phi-star-dv}). As revealed in Subsection \ref{X=Y=Z=V}, all these adjunctions are reduced to the ones described in \cite{Shen2013a,Shen2014} when $X=Y=Z=\sV$, which are now unified in the general framework of two-variable adjunction.

The main results of this paper are concerned about the representation theorems for fixed points of Isbell adjunctions and Kan adjunctions, whose prototypes can be traced in the fields of \emph{formal concept analysis} (FCA) \cite{Ganter1999,Davey2002} and \emph{rough set theory} (RST) \cite{Pawlak1982,Polkowski2002,Yao2004}. Explicitly, if $(X,Y,Z,\with,\lda,\rda)$ is a two-variable adjunction between \emph{complete} $\sV$-categories, then fixed points of $\uphi\dv\dphi$, $\psi^*\dv\psi_*$ and $\zeta_{\dag}\dv\zeta^{\dag}$ constitute complete $\sV$-categories
$$\Mphi,\quad\Kpsi\quad\text{and}\quad\Kdzeta,$$
respectively, which provide a (very general) categorical interpretation of the notion of \emph{concept lattice} in FCA and RST:
\begin{enumerate}[label=(\arabic*)]
\item If $X=Y=Z=\sV$, then $\Mphi$, $\Kpsi$ and $\Kdzeta$ are concept lattices valued in $\sV$ \cite{Bvelohlavek2001,Bvelohlavek2004,Lai2009,Shen2013a,Shen2014,Lai2017}; see Subsection \ref{X=Y=Z=V}.
\item If $\sV={\bf 2}$, the two-element Boolean algebra, then $\Mphi$, $\Kpsi$ and $\Kdzeta$ are \emph{multi-adjoint concept lattices} \cite{Medina2009,Medina2012,Lai2021}; see Subsection \ref{Multi-adjoint}.
\end{enumerate}
In particular, the representation theorems for $\Mphi$, $\Kpsi$ and $\Kdzeta$ in the above special cases have been established in the cited references. In our general setting, the desired representation theorems are presented as Theorem \ref{Mphi-representation}, Corollaries \ref{Kpsi-representation} and \ref{Kdzeta-representation}, which precisely characterize $\sV$-categories that are equivalent to $\Mphi$, $\Kpsi$ and $\Kdzeta$, respectively. It is worth pointing out that Theorem \ref{Mphi-representation} is not a straightforward generalization of \cite[Theorem 4.16]{Shen2013a} or \cite[Corollary 7.9]{Lai2017} for the case of $X=Y=Z=\sV$, since the absence of the Yoneda embedding in general $\sV$-categories of $\sV$-functors forces us to develop new techniques to complete the proof. The key tools are the $\sV$-bifunctors
$$\iota\colon A\ts X\to X^{A^{\op}}\quad\text{and}\quad\iotad\colon A\ts X^{\op}\to(X^A)^{\op}$$
constructed for any $\sV$-category $A$ and complete $\sV$-category $X$, which extend the notion of \emph{formal ball} in $\sV$-categories (Remark \ref{formal-ball}) and satisfy a generalized version of the Yoneda lemma (Propositions \ref{iota-dense} and \ref{gen-Yoneda}).

\section{$\sV$-categories} \label{V-categories}

Throughout, let
$$\sV=(\sV,\ts,k)$$
denote a \emph{commutative quantale} \cite{Rosenthal1990}; that is, $(\sV,\ts,k)$ is a commutative monoid with the underlying set $\sV$ being a complete lattice, such that
$$x\ts\Big(\bv_{i\in I}y_i\Big)=\bv_{i\in I}x\ts y_i$$
for all $x,y_i\in\sV$ $(i\in I)$. The right adjoint induced by the multiplication, denoted by $\ra$,
$$(x\ts -)\dv(x\ra -)\colon\sV\to\sV,$$
satisfies
$$x\ts y\leq z\iff x\leq y\ra z$$
for all $x,y,z\in\sV$.

A \emph{$\sV$-category} \cite{Lawvere1973,Kelly1982} consists of a set $X$ and a map $\al\colon X\times X\to\sV$ satisfying
$$k\leq\al(x,x)\quad\text{and}\quad\al(x,y)\ts\al(y,z)\leq\al(x,z)$$
for all $x,y,z\in X$. For simplicity, we abbreviate the pair $(X,\al)$ to $X$ and write $X(x,y)$ instead of $\al(x,y)$ if no confusion arises.

Every $\sV$-category $X$ is equipped with an underlying (pre)order given by
$$x\leq y\iff k\leq X(x,y)$$
for all $x,y\in X$. We write $x\cong y$ iff $x\leq y$ and $y\leq x$. A $\sV$-category $X$ is \emph{separated} (also \emph{skeletal}) if $x=y$ whenever $x\cong y$ in $X$.

A map $f\colon X\to Y$ between $\sV$-categories is a \emph{$\sV$-functor} if
$$X(x,y)\leq Y(fx,fy)$$
for all $x,y\in X$. With the order of $\sV$-functors given by
$$f\leq g\colon X\to Y\iff\forall x\in X\colon fx\leq gx\iff\forall x\in X\colon k\leq Y(fx,gx),$$
we obtain a 2-category
$$\VCat$$
of $\sV$-categories and $\sV$-functors.

Given a $\sV$-functor $f\colon X\to Y$, we say that:
\begin{itemize}
\item $f$ is \emph{fully faithful}, if $X(x,y)=Y(fx,fy)$ for all $x,y\in X$.
\item $f$ is \emph{essentially surjective}, if for each $y\in Y$ there exists $x\in X$ such that $y\cong fx$.
\item $f$ is an \emph{equivalence} (resp. \emph{isomorphism}) of $\sV$-categories, if there exists a $\sV$-functor $g\colon Y\to X$ with $gf\cong 1_X$ (resp. $gf=1_X$) and $fg\cong 1_Y$ (resp. $fg=1_Y$), where $1_X$ and $1_Y$ refer to the identity maps on $X$ and $Y$, respectively. In this case, $X$ and $Y$ are \emph{equivalent} (resp. \emph{isomorphic}) $\sV$-categories, and denoted by $X\simeq Y$ (resp. $X\cong Y$).
\end{itemize}

Assuming the axiom of choice, the following characterization of equivalences of $\sV$-categories is well known:

\begin{prop} \cite{Kelly1982} \label{f-equiv-ff-sur}
A $\sV$-functor $f\colon X\to Y$ is an equivalence of $\sV$-categories if, and only if, $f$ is fully faithful and essentially surjective.
\end{prop}

\begin{exmp} \label{V-cat-exmp}
Some basic examples of $\sV$-categories are listed below:
\begin{enumerate}[label=(\arabic*)]
\item \label{V-cat-exmp:V} $\sV$ itself is a separated $\sV$-category with
    $$\sV(x,y)=x\ra y$$ for all $x,y\in\sV$.
\item \label{V-cat-exmp:metric} Let $[0,\infty]=([0,\infty],\geq,+,0)$ denote Lawvere's quantale. Then $[0,\infty]$-categories are (generalized) metric spaces \cite{Lawvere1973}, and $[0,\infty]$-functors are \emph{non-expansive} maps.
\item \label{V-cat-exmp:dual} Every $\sV$-category $X$ has a \emph{dual} $X^{\op}$, which has the same underlying set as $X$, and
    $$X^{\op}(x,y)=X(y,x)$$
    for all $x,y\in X$.
\item \label{V-cat-exmp:product} The \emph{product} of $\sV$-categories $X$, $Y$, denoted by $X\times Y$, has the cartesian product of their underlying sets as its set of objects, and
    $$(X\times Y)((x,y),(x',y'))=X(x,x')\wedge Y(y,y')$$
    for all $x,x'\in X$, $y,y'\in Y$.
\item \label{V-cat-exmp:ts} The \emph{tensor product} of $\sV$-categories $X$, $Y$, denoted by $X\ts Y$, has the cartesian product of their underlying sets as its set of objects, and
    $$(X\ts Y)((x,y),(x',y'))=X(x,x')\ts Y(y,y')$$
    for all $x,x'\in X$, $y,y'\in Y$.
\item \label{V-cat-exmp:functor} Given $\sV$-categories $X$, $Y$, we denote by $Y^X$ the $\sV$-category of $\sV$-functors from $X$ to $Y$, with
    \begin{equation} \label{Y-X-f-g}
    Y^X(f,g)=\bw_{x\in X}Y(fx,gx)
    \end{equation}
    for all $f,g\in Y^X$. In particular, $\sV^{X^{\op}}$ and $(\sV^X)^{\op}$ are called the \emph{presheaf} and \emph{copresheaf $\sV$-categories} of $X$, respectively, with
    $$\sV^{X^{\op}}(\mu,\mu')=\bw_{x\in X}\mu x\ra\mu' x\quad\text{and}\quad(\sV^X)^{\op}(\lam,\lam')=\bw_{x\in X}\lam' x\ra\lam x$$
    for all $\mu,\mu'\in\sV^{X^{\op}}$, $\lam,\lam'\in\sV^X$.
\end{enumerate}
\end{exmp}

\begin{rem} \label{dual}
The dual of a $\sV$-category given in Example \ref{V-cat-exmp}\ref{V-cat-exmp:dual} may be expanded to an isomorphism
$$(-)^{\op}\colon\VCat\to(\VCat)^{\co}$$
of 2-categories, where ``$\co$'' refers to the dualization of 2-cells. Explicitly, the dual of a $\sV$-functor $f\colon X\to Y$, denoted by $f^{\op}\colon X^{\op}\to Y^{\op}$, is the same map as $f$ on objects, but $(f')^{\op}\leq f^{\op}$ whenever $f\leq f'\colon X\to Y$. Moreover, it is easy to see that
\begin{equation} \label{ts-dual}
(Y^X)^{\op}\cong(Y^{\op})^{X^{\op}}\quad\text{and}\quad(X\ts Y)^{\op}=X^{\op}\ts Y^{\op}
\end{equation}
for all $\sV$-categories $X$, $Y$.
\end{rem}

\section{Two-variable adjunctions in $\sV$-categories} \label{Two-variable-adjunction}

Recall that a pair of $\sV$-functors, $f\colon X\to Y$ and $g\colon Y\to X$, form an \emph{adjunction} in $\VCat$, denoted by $f\dv g$, if
$$1_X\leq gf\quad\text{and}\quad fg\leq 1_Y;$$
or equivalently, if
\begin{equation} \label{X-fx-y=Y-x-gy}
Y(fx,y)=X(x,gy)
\end{equation}
for all $x\in X$, $y\in Y$. It is well known that the $\sV$-functoriality of $f$, $g$ is implied by \eqref{X-fx-y=Y-x-gy}; that is, maps $f$, $g$ satisfying \eqref{X-fx-y=Y-x-gy} are necessarily $\sV$-functors (see, e.g., \cite[Theorem 2.10]{Lai2007}). In this case, $f$ is called a \emph{left adjoint} of $g$, and $g$ a \emph{right adjoint} of $f$.

The aim of this section is to introduce the notion of \emph{two-variable adjunction} in the context of $\sV$-categories, which appears in \cite{Gray1980,Riehl2016} for ordinary categories and in \cite{Shulman2006} for enriched categories. First, it is straightforward to check the following lemma, where a \emph{$\sV$-bifunctor} simply means a $\sV$-functor $X\ts Y\to Z$:

\begin{lem} \label{bifunctor}
For $\sV$-categories $X$, $Y$, $Z$, a map $\phi\colon X\ts Y\to Z$ is a $\sV$-bifunctor if, and only if,
\begin{enumerate}[label={\rm(\arabic*)}]
\item \label{bifunctor:x} $\phi(x,-)\colon Y\to Z$ is a $\sV$-functor for all $x\in X$, and
\item \label{bifunctor:y} $\phi(-,y)\colon X\to Z$ is a $\sV$-functor for all $y\in Y$.
\end{enumerate}
\end{lem}

%
%

Suppose that $X$, $Y$, $Z$ are $\sV$-categories and $\with\colon X\ts Y\to Z$ is a $\sV$-bifunctor. If the $\sV$-functors
$$-\with y\colon X\to Z\quad\text{and}\quad x\with -\colon Y\to Z$$
admit right adjoints in $\VCat$ for any $y\in Y$, $x\in X$, say,
\begin{equation} \label{lda-rda-def}
(-\with y)\dv(-\lda y)\colon Z\to X\quad\text{and}\quad (x\with -)\dv(x\rda -)\colon Z\to Y,
\end{equation}
then
\begin{equation} \label{with-lda-rda-adj}
Z(x\with y,z)=X(x,z\lda y)=Y(y,x\rda z)
\end{equation}
for all $x\in X$, $y\in Y$, $z\in Z$.

\begin{lem} \label{lda-rda-bifunctor}
$\lda\colon Z\ts Y^{\op}\to X$ and $\rda\colon X^{\op}\ts Z\to Y$ defined by \eqref{lda-rda-def} (or equivalently, by \eqref{with-lda-rda-adj}) are $\sV$-bifunctors.
\end{lem}

\begin{proof}
We check the $\sV$-bifunctoriality of $\lda$ as an example. By Lemma \ref{bifunctor}, it suffices to show that for each $z\in Z$,
$$z\lda -\colon Y^{\op}\to X$$
is a $\sV$-functor, which follows from
\begin{align*}
Y(y,y')&\leq Z((z\lda y')\with y,(z\lda y')\with y')\\
&\leq Z((z\lda y')\with y,(z\lda y')\with y')\otimes X(z\lda y',z\lda y')\\
&=Z((z\lda y')\with y,(z\lda y')\with y')\ts Z((z\lda y')\with y',z)\\
&\leq Z((z\lda y')\with y,z)\\
&=X(z\lda y',z\lda y),
\end{align*}
for all $y,y'\in Y$.
\end{proof}

\begin{defn} \label{two-var-adj-def} \cite{Shulman2006}
Let $X$, $Y$, $Z$ be $\sV$-categories. A \emph{two-variable adjunction} in $\VCat$ consists of $\sV$-bifunctors
$$\with\colon X\ts Y\to Z,\quad\lda\colon Z\ts Y^{\op}\to X,\quad\rda\colon X^{\op}\ts Z\to Y$$
such that
$$Z(x\with y,z)=X(x,z\lda y)=Y(y,x\rda z)$$
for all $x\in X$, $y\in Y$, $z\in Z$.
\end{defn}

For each $\sV$-category $X$, the \emph{Yoneda embedding} (resp. \emph{co-Yoneda embedding}) refers to the $\sV$-functor
$$\sy\colon X\to\sV^{X^{\op}},\quad x\mapsto X(-,x)\quad(\text{resp.}\ \syd\colon X\to(\sV^X)^{\op},\quad x\mapsto X(x,-)).$$
The following \emph{Yoneda lemma} is well known, which in particular implies that $\sy$ and $\syd$ are both fully faithful; however, $\sy$ and $\syd$ are injective maps only if $X$ is a separated $\sV$-category.

\begin{lem}[Yoneda] \cite{Kelly1982} \label{Yoneda}
For any $x\in X$, $\mu\in\sV^{X^{\op}}$, $\lam\in\sV^X$, it holds that
\begin{equation} \label{Yoneda-def}
\sV^{X^{\op}}(\sy x,\mu)=\mu x\quad\text{and}\quad(\sV^X)^{\op}(\lam,\syd x)=\lam x.
\end{equation}
\end{lem}

The \emph{supremum} (resp. \emph{infimum}) of $\mu\in\sV^{X^{\op}}$ (resp. $\lam\in\sV^X$), when it exists, is an object $\sup\mu\in X$ (resp. $\inf\lam\in X$) satisfying
\begin{align}
X(\sup\mu,x)&=\sV^{X^{\op}}(\mu,\sy x)=\bw_{y\in X}\mu y\ra X(y,x) \label{sup-def}\\
\Big(\text{resp.}\ \ X(x,\inf\lam)&=(\sV^X)^{\op}(\syd x,\lam)=\bw_{y\in X}\lam y\ra X(x,y)\Big). \label{inf-def}
\end{align}
for all $x\in X$. A $\sV$-category $X$ is \emph{cocomplete} if every $\mu\in\sV^{X^{\op}}$ admits a supremum, or equivalently, if $\sy\colon X\to\sV^{X^{\op}}$ has a left adjoint, given by $\sup\colon\sV^{X^{\op}}\to X$. It is well known that $X$ is cocomplete if and only if $X$ is \emph{complete} \cite{Stubbe2005}, where the completeness of $X$ is defined as $\syd\colon X\to(\sV^X)^{\op}$ admitting a right adjoint $\inf\colon(\sV^X)^{\op}\to X$.


Every $\sV$-functor $f\colon X\to Y$ induces a $\sV$-functor $f^{\ra}\colon\sV^X\to\sV^Y$ given by
\begin{equation} \label{fra-def}
(f^{\ra}\mu)y=\bv_{x\in X}Y(fx,y)\ts\mu x,
\end{equation}
and thus gives rise to $\sV$-functors
$$(f^{\op})^{\ra}\colon\sV^{X^{\op}}\to\sV^{Y^{\op}}\quad\text{and}\quad (f^{\ra})^{\op}\colon(\sV^X)^{\op}\to(\sV^Y)^{\op}$$
between (co)presheaf $\sV$-categories of $X$, $Y$.


\begin{prop} \cite{Stubbe2005} \label{left-adjoint-sup}
Let $f\colon X\to Y$ be a $\sV$-functor between complete $\sV$-categories. Then $f$ is a left (resp. right) adjoint in $\VCat$ if, and only if, $f$ preserves suprema (resp. infima) in the sense that
$$f{\sup}_X={\sup}_Y (f^{\op})^{\ra}\quad(\text{resp.}\ f{\inf}_X={\inf}_Y (f^{\ra})^{\op}).$$
\end{prop}

Therefore, if $X$, $Y$, $Z$ are complete $\sV$-categories, then a two-variable adjunction
$$(X,Y,Z,\with,\lda,\rda)$$
in $\VCat$ is completely determined by a $\sV$-bifunctor
$$\with\colon X\ts Y\to Z$$
that preserves suprema on both sides.

\section{Tensors and cotensors in $\sV$-categories} \label{Tensor}

In a $\sV$-category $X$, the \emph{tensor} (resp. \emph{cotensor}) \cite{Kelly1982} of $v\in\sV$ and $x\in X$, when it exists, is an object $v\star x\in X$ (resp. $v\rat x\in X$) such that
\begin{equation} \label{tensor-def}
X(v\star x,y)=v\ra X(x,y)\quad(\text{resp.}\ X(y,v\rat x)=v\ra X(y,x))
\end{equation}
for all $y\in X$. $X$ is \emph{tensored} (resp. \emph{cotensored}) if $v\star x$ (resp. $v\rat x$) exists for all $v\in\sV$, $x\in X$.

A $\sV$-category $X$ is \emph{order-complete} if the underlying ordered set of $X$ is complete.

\begin{prop} \cite{Stubbe2006} \label{V-complete-tensor}
A $\sV$-category $X$ is complete if, and only if, $X$ is tensored, cotensored and order-complete.
\end{prop}

\begin{prop} \cite{Stubbe2006} \label{left-adjoint-tensor}
Let $f\colon X\to Y$ be a $\sV$-functor between complete $\sV$-categories. Then $f$ is a left (resp. right) adjoint in $\VCat$ if, and only if, $f$ is a left (resp. right) adjoint between the underlying ordered sets of $X$, $Y$, and preserves tensors (resp. cotensors) in the sense that $f(v\star x)= v\star fx$ (resp. $f(v\rat x)=v\rat fx$) for all $v\in\sV$, $x\in X$.
\end{prop}

In fact, in a complete $\sV$-category $X$, \eqref{tensor-def} indicates that there are adjoint $\sV$-functors
\begin{equation} \label{tensor-adjoint}
(-\star x)\dv X(x,-)\colon X\to\sV\quad\text{and}\quad X(-,x)\dv(-\rat x)\colon\sV^{\op}\to X
\end{equation}
for all $x\in X$. Then it follows immediately from Proposition \ref{left-adjoint-tensor} that
\begin{equation} \label{X-x-bw-yi}
X\Big(x,\bw_{i\in I}y_i\Big)=\bw_{i\in I}X(x,y_i)\quad\text{and}\quad X\Big(\bv_{i\in I}y_i,x\Big)=\bw_{i\in I}X(y_i,x)
\end{equation}
for all $x,y_i\in X$ ($i\in I$), where $\displaystyle\bw\limits_{i\in I}y_i$ and $\displaystyle\bv\limits_{i\in I}y_i$ are calculated in the underlying order of $X$.

Let $X$ be a complete $\sV$-category. For each $\sV$-category $A$, we define
\begin{align}
&\iota(a,x)\colon A^{\op}\to X,\quad \iota(a,x)b=A(b,a)\star x, \label{iota-def}\\
&\iotad(a,x)\colon A\to X,\quad \iotad(a,x)b=A(a,b)\star x \label{iotad-def}
\end{align}
for all $a\in A$, $x\in X$.

\begin{lem} \label{iota-functor}
$\iota(a,x)\colon A^{\op}\to X$ and $\iotad(a,x)\colon A\to X$ are both $\sV$-functors.
\end{lem}

\begin{proof}
The $\sV$-functoriality of $\iota(a,x)$ follows from
\begin{align*}
A(b,c)&\leq A(c,a)\ra A(b,a)\\
&\leq\bw_{y\in Y}(A(b,a)\ra X(x,y))\ra(A(c,a)\ra X(x,y))\\
&=\bw_{y\in Y}X(A(b,a)\star x,y)\ra X(A(c,a)\star x,y)\\
&=X(A(c,a)\star x,\ A(b,a)\star x)\\
&=X(\iota(a,x)c,\iota(a,x)b)
\end{align*}
for all $b,c\in A$, and the $\sV$-functoriality of $\iotad(a,x)$ is obtained dually.
%
\end{proof}

\begin{rem} \label{formal-ball}
The $\sV$-functor $\iota(a,x)\colon A^{\op}\to X$ extends the notion of \emph{formal ball} in $\sV$-categories \cite{Kostanek2011}, which originates from (generalized) metric spaces (i.e., when $\sV=[0,\infty]$) \cite{Edalat1998,Rutten1998}. In fact, if $X=\sV$, then tensors in $\sV$ are given by
$$v\star r=v\ts r$$
for all $v,r\in\sV$. Consequently, for any $a\in A$ and $r\in\sV$, $\iota(a,r)\colon A^{\op}\to\sV$ is precisely a \emph{formal ball of center $a$ and radius $r$} in the sense of \cite[Definition 5.1]{Kostanek2011}.
\end{rem}

A $\sV$-functor $f\colon X\to Y$ is \emph{dense} (resp. \emph{codense}) if, for every $y\in Y$, there exists $\mu\in\sV^{X^{\op}}$ (resp. $\lam\in\sV^X$) such that
$$y=\sup (f^{\op})^{\ra}\mu\quad(\text{resp.}\ y=\inf(f^{\ra})^{\op}\lam).$$
The following properties of (co)dense $\sV$-functors are useful later: 

\begin{lem} (cf. \cite[Theorem 5.1]{Kelly1982} and \cite[Proposition 4.12]{Lai2017}) \label{dense-graph}
A $\sV$-functor $f\colon X\to Y$ is dense (resp. codense) if, and only if,
$$Y(y,y')=\bw_{x\in X}Y(fx,y)\ra Y(fx,y')\quad\Big(\text{resp.}\ Y(y,y')=\bw_{x\in X}Y(y',fx)\ra Y(y,fx)\Big)$$
for all $y,y'\in X$.
\end{lem}

\begin{lem} (cf. \cite[Proposition 5.9]{Kelly1982} and \cite[Corollary 4.13]{Lai2017}) \label{dense-comp}
\begin{enumerate}[label={\rm(\arabic*)}]
\item If a $\sV$-functor $f\colon X\to Y$ is essentially surjective, then $f$ is both dense and codense.
\item If $\sV$-functors $f\colon X\to Y$ and $g\colon Y\to Z$ are both dense (resp. codense), and $g$ is a left (resp. right) adjoint in $\VCat$, then $gf\colon X\to Z$ is dense (resp. codense).
\end{enumerate}
\end{lem}

Given a $\sV$-category $A$ and a complete $\sV$-category $X$, Lemma \ref{iota-functor} actually gives rise to well-defined maps
$$\iota\colon A\ts X\to X^{A^{\op}}\quad\text{and}\quad\iotad\colon A\ts X^{\op}\to(X^A)^{\op},$$
and moreover:

\begin{prop} \label{iota-dense}
$\iota\colon A\ts X\to X^{A^{\op}}$ is a dense $\sV$-bifunctor, and $\iotad\colon A\ts X^{\op}\to(X^A)^{\op}$ is a codense $\sV$-bifunctor.
\end{prop}

\begin{proof}
We only prove the claim for $\iota$, as the claim for $\iotad$ can be obtained dually.

First, $\iota$ is a $\sV$-bifunctor. This is because
\begin{align*}
(A\ts X)((a,x),(a',x'))&=A(a,a')\ts X(x,x')\\
&\leq(A(b,a)\ra A(b,a'))\ts X(x,x')\\
&\leq\bw_{b\in A}\bw_{y\in Y}(A(b,a')\ra X(x',y))\ra(A(b,a)\ra X(x,y))\\
&=\bw_{b\in A}\bw_{y\in Y}X(A(b,a')\star x',y)\ra X(A(b,a)\star x,y)\\
&=\bw_{b\in A}X(A(b,a)\star x,\ A(b,a')\star x')\\
&=X^{A^{\op}}(\iota(a,x),\iota(a',x'))
\end{align*}
for all $a,a'\in A$, $x,x'\in X$.

Second, $\iota$ is dense. By Lemma \ref{dense-graph} it suffices to show that
$$X^{A^{\op}}(\mu,\mu')=\bw_{a\in A}\bw_{x\in X}X^{A^{\op}}(\iota(a,x),\mu)\ra X^{A^{\op}}(\iota(a,x),\mu')$$
for all $\mu,\mu'\in X^{A^{\op}}$. To this end, note that for any $c\in A$, the $\sV$-functoriality of $\mu\colon A^{\op}\to X$ implies that
$$A(b,c)\leq X(\mu c,\mu b)$$
for all $b\in A$, and consequently
$$k\leq\bw_{b\in A}A(b,c)\ra X(\mu c,\mu b)=\bw_{b\in A}X(A(b,c)\star\mu c,\mu b)=X^{A^{\op}}(\iota(c,\mu c),\mu).$$
It follows that
\begin{align*}
&\bw_{a\in A}\bw_{x\in X}X^{A^{\op}}(\iota(a,x),\mu)\ra X^{A^{\op}}(\iota(a,x),\mu')\\
\leq{}&X^{A^{\op}}(\iota(c,\mu c),\mu)\ra X^{A^{\op}}(\iota(c,\mu c),\mu')\leq X^{A^{\op}}(\iota(c,\mu c),\mu')\\
\leq{}&X(\iota(c,\mu c)c,\mu'c)=X(A(c,c)\star\mu c,\mu'c)=A(c,c)\ra X(\mu c,\mu'c)\leq X(\mu c,\mu'c).
\end{align*}
Hence
$$\bw_{a\in A}\bw_{x\in X}X^{A^{\op}}(\iota(a,x),\mu)\ra X^{A^{\op}}(\iota(a,x),\mu')\leq\bw_{c\in A}X(\mu c,\mu'c)=X^{A^{\op}}(\mu,\mu'),$$
which is in fact an equation since the reverse inequality is trivial.
\end{proof}


\begin{prop} \label{gen-Yoneda}
Given a $\sV$-category $A$ and a complete $\sV$-category $X$, it holds that
\begin{equation} \label{gen-Yoneda-def}
X^{A^{\op}}(\iota(a,x),\mu)=X(x,\mu a)\quad\text{and}\quad(X^A)^{\op}(\lam,\iotad(a,x))=X(x,\lam a)
\end{equation}
for all $a\in A$, $x\in X$, $\mu\in X^{A^{\op}}$, $\lam\in X^A$.
\end{prop}

\begin{proof}
The identity for $\iota$ follows from
$$X^{A^{\op}}(\iota(a,x),\mu)=\bw_{b\in A}X(A(b,a)\star x,\mu b)=\bw_{b\in A}A(b,a)\ra X(x,\mu b)=X(x,\mu a)$$
for all $a\in A$, $x\in X$, $\mu\in X^{A^{\op}}$, and the identity for $\iotad$ follows dually.
\end{proof}

\begin{rem}
Proposition \ref{gen-Yoneda} may be viewed as a generalized version of the Yoneda lemma. In fact, if $X=\sV$, then for each $\sV$-category $A$, by setting
$$x=k$$
in \eqref{iota-def} and \eqref{iotad-def} we obtain the Yoneda embedding and the co-Yoneda embedding, i.e.,
$$\iota(a,k)=\sy_A a=A(-,a)\quad\text{and}\quad\iotad(a,k)=\syd_A a=A(a,-)$$
for all $a\in A$, and \eqref{gen-Yoneda-def} becomes exactly \eqref{Yoneda-def} in Lemma \ref{Yoneda}.
\end{rem}

\section{Isbell adjunctions and Kan adjunctions} \label{Isbell-Kan}

Let
$$(X,Y,Z,\with,\lda,\rda)$$
be a two-variable adjunction in $\VCat$. For each $\sV$-bifunctor
$$\phi\colon A^{\op}\ts B\to Z,$$
we define
\begin{align}
&\uphi\colon Y^{A^{\op}}\to(X^B)^{\op},\quad(\uphi\mu)b=\bw_{a\in A}\phi(a,b)\lda\mu a, \label{uphi-def}\\
&\dphi\colon (X^B)^{\op}\to Y^{A^{\op}},\quad(\dphi\lam)a=\bw_{b\in B}\lam b\rda\phi(a,b). \label{dphi-def}
\end{align}

\begin{prop} \label{uphi-dv-dphi}
Let $\phi\colon A^{\op}\ts B\to Z$ be a $\sV$-bifunctor. Then
$$\uphi\dv\dphi\colon (X^B)^{\op}\to Y^{A^{\op}}$$
in $\VCat$.
\end{prop}

\begin{proof}
First, $\uphi$ and $\dphi$ are well-defined maps. Given $\mu\in Y^{A^{\op}}$, the $\sV$-functoriality of $\uphi\mu\colon B\to X$ follows from
\begin{align*}
B(b,b')&\leq\bw_{a\in A}Z(\phi(a,b),\phi(a,b'))&(\text{Lemma \ref{bifunctor}})\\
&\leq\bw_{a\in A}X(\phi(a,b)\lda\mu a,\ \phi(a,b')\lda\mu a)&(\text{Lemma \ref{lda-rda-bifunctor}})\\
&\leq\bw_{a\in A}X\Big(\bw_{a'\in A}\phi(a',b)\lda\mu a',\ \phi(a,b')\lda\mu a\Big)&(\text{Equations \eqref{X-x-bw-yi}})\\
&=X\Big(\bw_{a'\in A}\phi(a',b)\lda\mu a',\bw_{a\in A}\phi(a,b')\lda\mu a\Big)&(\text{Equations \eqref{X-x-bw-yi}})\\
&=X((\uphi\mu)b,(\uphi\mu)b')&(\text{Equation \eqref{uphi-def}})
\end{align*}
for all $b,b'\in B$; and similarly, $\dphi\lam\colon A^{\op}\to Y$ is a $\sV$-functor whenever $\lam\in X^B$.


Second, $\uphi$ and $\dphi$ are $\sV$-functors and $\uphi\dv\dphi$ in $\VCat$. Indeed,
\begin{align*}
(X^B)^{\op}(\uphi\mu,\lam)&=\bw_{b\in B}X(\lam b,(\uphi\mu)b)\\
&=\bw_{b\in B}X\Big(\lam b,\bw_{a\in A}\phi(a,b)\lda\mu a\Big)&(\text{Equation \eqref{uphi-def}})\\
&=\bw_{b\in B}\bw_{a\in A}X(\lam b,\phi(a,b)\lda\mu a)&(\text{Equations \eqref{X-x-bw-yi}})\\
&=\bw_{a\in A}\bw_{b\in B}Y(\mu a,\lam b\rda\phi(a,b))&(\text{Definition \ref{two-var-adj-def}})\\
&=\bw_{a\in A}Y\Big(\mu a,\bw_{b\in B}\lam b\rda\phi(a,b)\Big)&(\text{Equations \eqref{X-x-bw-yi}})\\
&=\bw_{a\in A}Y(\mu a,(\dphi\lam)a)&(\text{Equation \eqref{dphi-def}})\\
&=Y^{A^{\op}}(\mu,\dphi\lam)
\end{align*}
for all $\mu\in Y^{A^{\op}}$, $\lam\in X^B$, which completes the proof.
\end{proof}

For each $\sV$-bifunctor $\phi\colon X\ts Y\to Z$, let $\phi^{\partial}\colon Y\ts X\to Z$ denote the $\sV$-bifunctor given by
$$\phi^{\partial}(y,x)=\phi(x,y)$$
for all $x\in X$, $y\in Y$.

\begin{lem} \label{two-var-adj-dual}
For $\sV$-bifunctors $\with\colon X\ts Y\to Z$, $\lda\colon Z\ts Y^{\op}\to X$, $\rda\colon X^{\op}\ts Z\to Y$, the following statements are equivalent:
\begin{enumerate}[label={\rm(\roman*)}]
\item \label{two-var-adj-dual:o} $(X,Y,Z,\with,\lda,\rda)$ is a two-variable adjunction in $\VCat$.
\item \label{two-var-adj-dual:rda} $(X,Z^{\op},Y^{\op},\rda,\lda^{\partial},\with)$ is a two-variable adjunction in $\VCat$.
\item \label{two-var-adj-dual:lda} $(Z^{\op},Y,X^{\op},\lda,\with,\rda^{\partial})$ is a two-variable adjunction in $\VCat$.
\end{enumerate}
\end{lem}

\begin{proof}
This is an immediate consequence of
\begin{align*}
&Z(x\with y,z)=X(x,z\lda y)=Y(y,x\rda z)\\
\iff{}&Y^{\op}(x\rda z,y)=X(x,y\lda^{\partial}z)=Z^{\op}(z,x\with y)\\
\iff{}&X^{\op}(z\lda y,x)=Z^{\op}(z,x\with y)=Y(y,z\rda^{\partial}x)
\end{align*}
for all $x\in X$, $y\in Y$, $z\in Z$.
\end{proof}

Note that the duality of $\sV$-categories (cf. Remark \ref{dual}) allows us to reformulate Proposition \ref{uphi-dv-dphi} as follows: the dual of each $\sV$-bifunctor
$$\phi\colon A^{\op}\ts B\to Z^{\op},\quad\text{i.e.},\quad\phi^{\op}\colon (A^{\op})^{\op}\ts B^{\op}\to Z,$$
induces an adjunction
$$\uphiop\dv\dphiop\colon(X^{B^{\op}})^{\op}\to Y^{(A^{\op})^{\op}}$$
in $\VCat$, given by \eqref{uphi-def} and \eqref{dphi-def}, which by duality corresponds to
\begin{equation} \label{uphi-dv-dphi-dual}
(\dphiop)^{\op}\dv(\uphiop)^{\op}\colon(Y^A)^{\op}\to X^{B^{\op}}.
\end{equation}
Applying \eqref{uphi-dv-dphi-dual} to the two-variable adjunctions of Lemma \ref{two-var-adj-dual}\ref{two-var-adj-dual:rda}\ref{two-var-adj-dual:lda}, we obtain that any $\sV$-bifunctors
$$\psi\colon A^{\op}\ts B\to Y\quad\text{and}\quad\zeta\colon A^{\op}\ts B\to X$$
give rise to adjunctions given in the following Proposition \ref{phi-star-dv}, where
\begin{align}
&\psi^*\colon X^{B^{\op}}\to Z^{A^{\op}},\quad(\psi^*\lam)a=\bv_{b\in B}\lam b\with\psi(a,b), \label{psi-ustar-def}\\
&\psi_*\colon Z^{A^{\op}}\to X^{B^{\op}},\quad(\psi_*\mu)b=\bw_{a\in A}\mu a\lda\psi(a,b) \label{psi-lstar-def}
\end{align}
and
\begin{align}
&\zeta_{\dag}\colon (Z^B)^{\op}\to (Y^A)^{\op},\quad(\zeta_{\dag}\lam)a=\bw_{b\in B}\zeta(a,b)\rda\lam b, \label{zeta-ldag-def}\\
&\zeta^{\dag}\colon(Y^A)^{\op}\to(Z^B)^{\op},\quad(\zeta^{\dag}\mu)b=\bv_{a\in A}\zeta(a,b)\with\mu a. \label{zeta-udag-def}
\end{align}

\begin{prop} \label{phi-star-dv}
Let $\psi\colon A^{\op}\ts B\to Y$ and $\zeta\colon A^{\op}\ts B\to X$ be $\sV$-bifunctors. Then
$$\psi^*\dv\psi_*\colon Z^{A^{\op}}\to X^{B^{\op}}\quad\text{and}\quad\zeta_{\dag}\dv\zeta^{\dag}\colon (Y^A)^{\op}\to(Z^B)^{\op}$$
in $\VCat$.
\end{prop}


As we will see in Subsection \ref{X=Y=Z=V}, adjunctions given by Propositions \ref{uphi-dv-dphi} and \ref{phi-star-dv} generalize Isbell adjunctions and Kan adjunctions induced by \emph{$\sV$-distributors} between $\sV$-categories, respectively. Hence, following the terminologies of \cite{Shen2013a,Shen2014}, adjunctions of the forms $\uphi\dv\dphi$ and $\psi^*\dv\psi_*$, $\zeta_{\dag}\dv\zeta^{\dag}$ are called \emph{Isbell adjunctions} and \emph{Kan adjunctions}, respectively.

\section{Representation theorems} \label{Representation}

In this section, we assume that $X$, $Y$, $Z$ are \emph{complete} $\sV$-categories, and
$$(X,Y,Z,\with,\lda,\rda)$$
is a two-variable adjunction in $\VCat$.

Given a $\sV$-category $A$, a $\sV$-functor $h\colon A\to A$ is a \emph{$\sV$-closure operator} if
$$1_A\leq h\quad\text{and}\quad hh\cong h.$$
In particular, each adjunction $f\dv g\colon B\to A$ in $\VCat$ induces a $\sV$-closure operator $gf\colon A\to A$. We denote by
$$\Fix(h):=\{a\in A\mid ha\cong a\}$$
the $\sV$-subcategory of $A$ consisting of fixed points\footnote{Strictly speaking, ``fixed point'' should read ``pseudo-fixed point'' here, since $a\in\Fix(h)$ satisfies $ha\cong a$ instead of $ha=a$.} of $h$.

\begin{prop} \label{Fix-h-adjoint-complete} {\rm\cite{Shen2013a}}
Let $h\colon A\to A$ be a $\sV$-closure operator. Then
\begin{itemize}
\item[\rm (1)] the inclusion $\sV$-functor $\Fix(h)\ \to/^(->/A$ is the right adjoint of the codomain restriction $h\colon A\to\Fix(h)$;
\item[\rm (2)] $\Fix(h)$ is a complete $\sV$-category provided that $A$ is complete.
\end{itemize}
\end{prop}

If $X$ is a complete $\sV$-category, then for any $\sV$-category $A$, it is easy to see that $X^A$ is equipped with the pointwise tensors, cotensors and underlying joins inherited from $X$, and thus $X^A$ is also a complete $\sV$-category by Proposition \ref{V-complete-tensor}. Therefore, by Proposition \ref{uphi-dv-dphi}, every $\sV$-bifunctor
$$\phi\colon A^{\op}\ts B\to Z$$
gives rise to a complete $\sV$-category
$$\Mphi:=\Fix(\dphi\uphi)=\{\mu\in Y^{A^{\op}}\mid\dphi\uphi\mu\cong\mu\}$$
of the fixed points of the $\sV$-closure operator $\dphi\uphi\colon Y^{A^{\op}}\to Y^{A^{\op}}$ induced by the Isbell adjunction $\uphi\dv\dphi$.

The following theorem is the main result of this paper, which characterizes those $\sV$-categories that are equivalent to $\Mphi$:

\begin{thm} \label{Mphi-representation}
Let $\phi\colon A^{\op}\ts B\to Z$ be a $\sV$-bifunctor. A $\sV$-category $C$ is equivalent to $\Mphi$ if, and only if, $C$ is complete and there exist a dense $\sV$-bifunctor $\al\colon A\ts Y\to C$ and a codense $\sV$-bifunctor $\be\colon B\ts X^{\op}\to C$ such that
\begin{equation} \label{al-be-Z-phi}
C(\al(a,y),\be(b,x))=Z(x\with y,\phi(a,b))
\end{equation}
for all $a\in A$, $b\in B$, $x\in X$, $y\in Y$.
\end{thm}

Before proving this theorem, we present the following lemma as a preparation:

\begin{lem} \label{uphi-iota}
For any $a\in A$, $b\in B$, $x\in X$, $y\in Y$, it holds that
$$(\uphi\iota_{A,Y}(a,y))b\cong\phi(a,b)\lda y\quad\text{and}\quad(\dphi\iotad_{B,X}(b,x))a\cong x\rda\phi(a,b).$$
\end{lem}

\begin{proof}
The first isomorphism follows from
\begin{align*}
X(x',(\uphi\iota_{A,Y}(a,y))b)&=X\Big(x',\bw_{a'\in A}\phi(a',b)\lda(A(a',a)\star y)\Big)&(\text{Equations \eqref{iota-def} and \eqref{uphi-def}})\\
&=\bw_{a'\in A}X(x',\phi(a',b)\lda(A(a',a)\star y))&(\text{Equations \eqref{X-x-bw-yi}})\\
&=\bw_{a'\in A}Y(A(a',a)\star y,\ x'\rda\phi(a',b))&(\text{Definition \ref{two-var-adj-def}})\\
&=\bw_{a'\in A}A(a',a)\ra Y(y,x'\rda\phi(a',b))&(\text{Equations \eqref{tensor-def}})\\
&=Y(y,x'\rda\phi(a,b))\\
&=X(x',\phi(a,b)\lda y)&(\text{Definition \ref{two-var-adj-def}})
\end{align*}
for all $x'\in X$, where the penultimate equality holds since
$$A(a',a)\leq Z(\phi(a,b),\phi(a',b))\leq Y(x'\rda\phi(a,b),x'\rda\phi(a',b))$$
by applying Lemmas \ref{bifunctor} and \ref{lda-rda-bifunctor} to $\phi\colon A^{\op}\ts B\to Z$ and $\rda\colon X^{\op}\ts Z\to Y$. The second isomorphism is obtained dually.
%
\end{proof}

\begin{proof}[Proof of Theorem \ref{Mphi-representation}]
{\bf Necessity.} It suffices to prove the case $C=\Mphi$. We show that
$$\al\colon A\ts Y\to\Mphi\quad\text{and}\quad\be\colon B\ts X^{\op}\to\Mphi$$
given by
$$\al(a,y)=\dphi\uphi\iota_{A,Y}(a,y)\quad\text{and}\quad\be(b,x)=\dphi\iotad_{B,X}(b,x)$$
are dense and codense, respectively, and satisfy \eqref{al-be-Z-phi}.

First, the density of $\al$ and the codensity of $\be$ both follow from Lemma \ref{dense-comp}. Indeed, $\al$ is the composition of the dense $\sV$-bifunctor $\iota_{A,Y}\colon A\ts Y\to Y^{A^{\op}}$ (Proposition \ref{iota-dense}) and the codomain restriction of the $\sV$-closure operator $\dphi\uphi\colon Y^{A^{\op}}\to Y^{A^{\op}}$ (cf. Proposition \ref{Fix-h-adjoint-complete}), while $\be$ is the composition of the codense $\sV$-bifunctor $\iotad_{B,X}\colon B\ts X^{\op}\to(X^B)^{\op}$ (Proposition \ref{iota-dense}) and the codomain restriction of the right adjoint $\sV$-functor $\dphi\colon (X^B)^{\op}\to Y^{A^{\op}}$.

Second, Equation \eqref{al-be-Z-phi} follows from
\begin{align*}
\Mphi(\al(a,y),\be(b,x))&=Y^{A^{\op}}(\dphi\uphi\iota_{A,Y}(a,y),\dphi\iotad_{B,X}(b,x))\\
&=(X^B)^{\op}(\uphi\dphi\uphi\iota_{A,Y}(a,y),\iotad_{B,X}(b,x))&(\uphi\dv\dphi)\\
&=(X^B)^{\op}(\uphi\iota_{A,Y}(a,y),\iotad_{B,X}(b,x))&(\uphi\dv\dphi)\\
&=X(x,(\uphi\iota_{A,Y}(a,y))b)&(\text{Proposition \ref{gen-Yoneda}})\\
&=X(x,\phi(a,b)\lda y)&(\text{Lemma \ref{uphi-iota}})\\
&=Z(x\with y,\phi(a,b)) &(\text{Definition \ref{two-var-adj-def}})
\end{align*}
for all $a\in A$, $b\in B$, $x\in X$, $y\in Y$.

{\bf Sufficiency.} We show that
$$h\colon\Mphi\to C,\quad h\mu=\bv_{a\in A}\al(a,\mu a)$$
is a fully faithful and essentially surjective $\sV$-functor, and thus an equivalence of $\sV$-categories by Proposition \ref{f-equiv-ff-sur}.

First, $h$ is a fully faithful $\sV$-functor. Note that
\begin{equation} \label{h-mu-be-b-x}
C(h\mu,\be(b,x))=X(x,(\uphi\mu)b)
\end{equation}
for all $\mu\in\Mphi$, $b\in B$, $x\in X$, because
\begin{align*}
C(h\mu,\be(b,x))&=C\Big(\bv_{a\in A}\al(a,\mu a),\be(b,x)\Big)\\
&=\bw_{a\in A}C(\al(a,\mu a),\be(b,x))&(\text{Equations \eqref{X-x-bw-yi}})\\
&=\bw_{a\in A}Z(x\with\mu a,\phi(a,b))&(\text{Equation \eqref{al-be-Z-phi}})\\
&=\bw_{a\in A}X(x,\phi(a,b)\lda\mu a)&(\text{Definition \ref{two-var-adj-def}})\\
&=X\Big(x,\bw_{a\in A}\phi(a,b)\lda\mu a\Big)&(\text{Equations \eqref{X-x-bw-yi}})\\
&=X(x,(\uphi\mu)b).&(\text{Equation \eqref{uphi-def}})
\end{align*}
It follows that
\begin{align*}
\Mphi(\mu,\mu')&=Y^{A^{\op}}(\mu,\dphi\uphi\mu')&(\mu'\in\Mphi)\\
&=(X^B)^{\op}(\uphi\mu,\uphi\mu')&(\uphi\dv\dphi)\\
&=\bw_{b\in B}X((\uphi\mu')b,(\uphi\mu)b)\\
&=\bw_{b\in B}\bw_{x\in X}X(x,(\uphi\mu')b)\ra X(x,(\uphi\mu)b)\\
&=\bw_{b\in B}\bw_{x\in X}C(h\mu',\be(b,x))\ra C(h\mu,\be(b,x))&(\text{Equation \eqref{h-mu-be-b-x}})\\
&=C(h\mu,h\mu')&(\text{Lemma \ref{dense-graph}})
\end{align*}
for all $\mu,\mu'\in\Mphi$, where Lemma \ref{dense-graph} is applied to the codense $\sV$-bifunctor $\be\colon B\ts X^{\op}\to C$ in the last equality.

Second, $h$ is essentially surjective. For any $c\in C$, note that $C(\al(a,-),c)\in\sV^{Y^{\op}}$ for all $a\in A$. Thus it makes sense to define
$$\mu\colon A^{\op}\to Y,\quad \mu a={\sup}_Y C(\al(a,-),c),$$
which is a $\sV$-functor because
\begin{align*}
A(a,a')&\leq\bw_{y\in Y}C(\al(a,y),\al(a',y))&(\sV\text{-functoriality of}\ \al)\\
&\leq\bw_{y\in Y}C(\al(a',y),c)\ra C(\al(a,y),c)\\
&=V^{Y^{\op}}(C(\al(a',-),c),C(\al(a,-),c))\\
&\leq Y(\mu a',\mu a)&(\sV\text{-functoriality of}\ {\sup}_Y)
\end{align*}
for all $a,a'\in A$. We claim that
$$c\cong h\dphi\uphi\mu.$$
Indeed,
\begin{equation} \label{C-h-dphi-uphi-mu}
C(c,\be(b,x))=C(h\dphi\uphi\mu,\be(b,x))
\end{equation}
for all $b\in B$, $x\in X$, because the density of $\al\colon A\ts Y\to C$ guarantees that
\begin{align*}
C(c,\be(b,x))&=\bw_{a\in A}\bw_{y\in Y}C(\al(a,y),c)\ra C(\al(a,y),\be(b,x))&(\text{Lemma \ref{dense-graph}})\\
&=\bw_{a\in A}\bw_{y\in Y}C(\al(a,y),c)\ra Z(x\with y,\phi(a,b))&(\text{Equation \eqref{al-be-Z-phi}})\\
&=\bw_{a\in A}\bw_{y\in Y}C(\al(a,y),c)\ra Y(y,x\rda\phi(a,b))&(\text{Definition \ref{two-var-adj-def}})\\
&=\bw_{a\in A}Y(\mu a,x\rda\phi(a,b))&(\text{Equation \eqref{sup-def}})\\
&=\bw_{a\in A}X(x,\phi(a,b)\lda\mu a)&(\text{Definition \ref{two-var-adj-def}})\\
&=X\Big(x,\bw_{a\in A}\phi(a,b)\lda\mu a\Big)&(\text{Equations \eqref{X-x-bw-yi}})\\
&=X(x,(\uphi\mu)b)&(\text{Equation \eqref{uphi-def}})\\
&=X(x,(\uphi\dphi\uphi\mu)b)&(\uphi\dv\dphi)\\
&=C(h\dphi\uphi\mu,\be(b,x)).&(\text{Equation \eqref{h-mu-be-b-x}})
\end{align*}
Hence, by the codensity of $\be\colon B\ts X^{\op}\to C$,
\begin{align*}
C(c,c')&=\bw_{b\in B}\bw_{x\in X}C(c',\be(b,x))\ra C(c,\be(b,x))&(\text{Lemma \ref{dense-graph}})\\
&=\bw_{b\in B}\bw_{x\in X}C(c',\be(b,x))\ra C(h\dphi\uphi\mu,\be(b,x))&(\text{Equation \eqref{C-h-dphi-uphi-mu}})\\
&=C(h\dphi\uphi\mu,c')&(\text{Lemma \ref{dense-graph}})
\end{align*}
for all $c'\in C$, which completes the proof.
\end{proof}

\begin{rem}
The above proof for Theorem \ref{Mphi-representation} is a direct one. Besides, an easier but indirect approach to the sufficiency of the condition in Theorem \ref{Mphi-representation} may be formulated by applying the results of \cite{Lai2017}. Explicitly, since we have dense $\sV$-functors $\al$, $\iota_{A,Y}$ and codense $\sV$-functors $\be$, $\iotad_{B,X}$ satisfying
$$(X^B)^{\op}(\uphi\iota_{A,Y}(a,y),\iotad_{B,X}(b,x))=Z(x\with y,\phi(a,b))=C(\al(a,y),\be(b,x))$$
$$\bfig
\morphism(-1000,0)<600,300>[A\ts Y`Y^{A^{\op}};\iota_{A,Y}]
\morphism(1000,0)<-600,300>[B\ts X^{\op}`(X^B)^{\op};\iotad_{B,X}]
\morphism(-1000,0)|b|<1000,-250>[A\ts Y`C;\al]
\morphism(1000,0)|b|<-1000,-250>[B\ts X^{\op}`C;\be]
\morphism(-400,300)|a|/@{>}@<4pt>/<800,0>[Y^{A^{\op}}`(X^B)^{\op};\uphi]
\morphism(400,300)|b|/@{>}@<4pt>/<-800,0>[(X^B)^{\op}`Y^{A^{\op}};\dphi]
\place(-30,300)[\mbox{\scriptsize$\bot$}]
\efig$$
for all $a\in A$, $b\in B$, $x\in X$, $y\in Y$, it follows from \cite[Theorem 5.1]{Lai2017} that $C$ is equivalent to $\Mphi$.
\end{rem}

Now, in order to derive representation theorems for the complete $\sV$-categories of fixed points of Kan adjunctions given by Proposition \ref{phi-star-dv}, i.e.,
\begin{align*}
\Kpsi&:=\Fix(\psi_*\psi^*)=\{\lam\in X^{B^{\op}}\mid\psi_*\psi^*\lam\cong\lam\},\\
\Kdzeta&:=\Fix(\zeta^{\dag}\zeta_{\dag})=\{\lam\in Z^B\mid \zeta^{\dag}\zeta_{\dag}\lam\cong\lam\},
\end{align*}
let us apply Theorem \ref{Mphi-representation} to the adjunction \eqref{uphi-dv-dphi-dual}. Indeed, for any $\sV$-bifunctor $\phi\colon A^{\op}\ts B\to Z^{\op}$, the $\sV$-category of fixed points of \eqref{uphi-dv-dphi-dual} is precisely $(\Mphi^{\op})^{\op}$. By Theorem \ref{Mphi-representation}, the dual of a complete $\sV$-category $C$ is equivalent to $\Mphi^{\op}$ if, and only if, there exist a dense $\sV$-bifunctor $\al\colon A^{\op}\ts Y\to C^{\op}$ and a codense $\sV$-bifunctor $\be\colon B^{\op}\ts X^{\op}\to C^{\op}$ with
$$C(\be(b,x),\al(a,y))=C^{\op}(\al(a,y),\be(b,x))=Z(x\with y,\phi(a,b))$$
for all $a\in A$, $b\in B$, $x\in X$, $y\in Y$. Therefore, by the duality of $\sV$-categories we have:

\begin{cor} \label{Mphiop-representation}
Let $\phi\colon A^{\op}\ts B\to Z^{\op}$ be a $\sV$-bifunctor. A $\sV$-category $C$ is equivalent to $(\Mphi^{\op})^{\op}$ if, and only if, $C$ is complete and there exist a dense $\sV$-bifunctor $\be\colon B\ts X\to C$ and a codense $\sV$-bifunctor $\al\colon A\ts Y^{\op}\to C$ such that
$$C(\be(b,x),\al(a,y))=Z(x\with y,\phi(a,b))$$
for all $a\in A$, $b\in B$, $x\in X$, $y\in Y$.
\end{cor}

The representation theorems for $\Kpsi$ and $\Kdzeta$ are then obtained by applying Corollary \ref{Mphiop-representation} to Proposition \ref{phi-star-dv}:

\begin{cor} \label{Kpsi-representation}
Let $\psi\colon A^{\op}\ts B\to Y$ be a $\sV$-bifunctor. A $\sV$-category $C$ is equivalent to $\Kpsi$ if, and only if, $C$ is complete and there exist a dense $\sV$-bifunctor $\be\colon B\ts X\to C$ and a codense $\sV$-bifunctor $\al\colon A\ts Z\to C$ such that
$$C(\be(b,x),\al(a,z))=Y(\psi(a,b),x\rda z)$$
for all $a\in A$, $b\in B$, $x\in X$, $z\in Z$.
\end{cor}

\begin{cor} \label{Kdzeta-representation}
Let $\zeta\colon A^{\op}\ts B\to X$ be a $\sV$-bifunctor. A $\sV$-category $C$ is equivalent to $\Kdzeta$ if, and only if, $C$ is complete and there exist a dense $\sV$-bifunctor $\be\colon B\ts Z^{\op}\to C$ and a codense $\sV$-bifunctor $\al\colon A\ts Y^{\op}\to C$ such that
$$C(\be(b,z),\al(a,y))=X(\zeta(a,b),z\lda y)$$
for all $a\in A$, $b\in B$, $y\in Y$, $z\in Z$.
\end{cor}

\section{Examples}

\subsection{When $X=Y=Z=\sV$} \label{X=Y=Z=V}

When $X=Y=Z=\sV$, the multiplication $\ts$ of the quantale $\sV$ obviously induces a two-variable adjunction on $\sV$, since
$$V(x\ts y,z)=V(x,y\ra z)=V(y,x\ra z)$$
for all $x,y,z\in\sV$. A $\sV$-bifunctor $\phi\colon A^{\op}\ts B\to\sV$ is precisely a \emph{$\sV$-distributor}
$$\phi\colon A\oto B;$$
that is, a map $\phi\colon A\times B\to\sV$ satisfying
$$B(b,b')\ts\phi(a,b)\ts A(a',a)\leq\phi(a',b')$$
for all $a,a'\in A$, $b,b'\in B$. Therefore, the Isbell adjunction $\uphi\dv\dphi$ and Kan adjunctions $\phi^*\dv\phi_*$, $\phi^{\dag}\dv\phi_{\dag}$ induced by $\phi$ are exactly the ones given in \cite[Proposition 4.1]{Shen2013a} and \cite[Proposition 6.2.1]{Shen2014}, respectively, when $\CQ=\sV$ is a commutative quantale. In this case, $\Mphi$, $\Kphi$ and $\Kdphi$ are separated and complete $\sV$-categories since so is $\sV$, and the corresponding representation theorems are formulated as follows:

\begin{cor} \label{Mphi-Kphi-representation-V}
Let $\phi\colon A\oto B$ be a $\sV$-distributor, and let $C$ be a separated $\sV$-category.
\begin{enumerate}[label={\rm(\arabic*)}]
\item \label{Mphi-Kphi-representation-V:Mphi} $C$ is isomorphic to $\Mphi$ if, and only if, $C$ is complete and there exist a dense $\sV$-bifunctor $\al\colon A\ts\sV\to C$ and a codense $\sV$-bifunctor $\be\colon B\ts\sV^{\op}\to C$ and such that
    $$C(\al(a,y),\be(b,x))=(x\ts y)\ra\phi(a,b)$$
    for all $a\in A$, $b\in B$, $x,y\in\sV$.
\item \label{Mphi-Kphi-representation-V:Kphi} $C$ is isomorphic to $\Kphi$ if, and only if, $C$ is complete and there exist a dense $\sV$-bifunctor $\be\colon B\ts\sV\to C$ and a codense $\sV$-bifunctor $\al\colon A\ts\sV\to C$ such that
    $$C(\be(b,x),\al(a,z))=\phi(a,b)\ra(x\ra z)$$
    for all $a\in A$, $b\in B$, $x,z\in\sV$.
\item \label{Mphi-Kphi-representation-V:Kdphi} $C$ is isomorphic to $\Kdphi$ if, and only if, $C$ is complete and there exist a dense $\sV$-bifunctor $\be\colon B\ts\sV^{\op}\to C$ and a codense $\sV$-bifunctor $\al\colon A\ts\sV^{\op}\to C$ such that
    $$C(\be(b,z),\al(a,y))=\phi(a,b)\ra(y\ra z)$$
    for all $a\in A$, $b\in B$, $y,z\in\sV$.
\end{enumerate}
\end{cor}

\begin{rem}
As elaborated in \cite[Section 7]{Lai2017}, \ref{Mphi-Kphi-representation-V:Mphi} and \ref{Mphi-Kphi-representation-V:Kphi} of Corollary \ref{Mphi-Kphi-representation-V} are actually the commutative case of \cite[Corollary 7.9]{Lai2017}, whose prototypes come from \cite[Theorem 14(2)]{Bvelohlavek2004} and \cite[Proposition 7.3]{Popescu2004}, respectively.
\end{rem}

\subsection{Multi-adjoint concept lattices} \label{Multi-adjoint}

In the special case of $\sV={\bf 2}$, the two-element Boolean algebra, a two-variable adjunction with respect to ordered sets $X$, $Y$, $Z$ consists of maps
$$\with\colon X\times Y\to Z,\quad\lda\colon Z\times Y\to X,\quad\rda\colon X\times Z\to Y$$
such that
$$x\with y\leq z\iff x\leq z\lda y\iff y\leq x\rda z$$
for all $x\in X$, $y\in Y$, $z\in Z$; that is, $(X,Y,Z,\with,\lda,\rda)$ is an \emph{adjoint triple} in the sense of \cite{Medina2009,Medina2012}. For any sets $A$, $B$ and maps
$$\phi\colon A\times B\to Z,\quad\psi\colon A\times B\to Y\quad\text{and}\quad\zeta\colon A\times B\to X,$$
\begin{itemize}
\item $\Mphi$ is the \emph{multi-adjoint concept lattice} of the context $(A,B,\phi)$ \cite{Medina2009,Lai2021},
\item $\Kpsi$ is the \emph{multi-adjoint property-oriented concept lattice} of the context $(A,B,\psi)$ \cite{Medina2012,Lai2021},
\item $\Kdzeta$ is the \emph{multi-adjoint object-oriented concept lattice} of the context $(A,B,\zeta)$ \cite{Medina2012,Lai2021}.
\end{itemize}
In particular, our Theorem \ref{Mphi-representation} formalizes the representation theorem of the multi-adjoint concept lattice (see \cite[Theorem 20]{Medina2009}) when $\sV={\bf 2}$.

\subsection{Tensors and cotensors} \label{Tensor-Cotensor}

For every complete $\sV$-category $X$, from \eqref{tensor-adjoint} we see that its tensors and cotensors give rise to a two-variable adjunction
\begin{equation} \label{tensor-cotensor-two-variable}
(\sV,X,X,\star,X^{\op}(-,-),\rat);
\end{equation}
explicitly,
$$X(v\star x,y)=\sV(v,X^{\op}(y,x))=X(x,v\rat y)$$
for all $v\in\sV$, $x,y\in X$. In this case:
\begin{enumerate}[label=(\arabic*)]
\item The Isbell adjunction induced by a $\sV$-bifunctor $\phi\colon A^{\op}\ts B\to X$ is given by
    \begin{align}
    &\uphi\colon X^{A^{\op}}\to(\sV^B)^{\op},\quad(\uphi\mu)b=\bw_{a\in A}X(\mu a,\phi(a,b))=X^{A^{\op}}(\mu,\phi(-,b)),\label{uphi-tensor-cotensor-def}\\
    &\dphi\colon (\sV^B)^{\op}\to X^{A^{\op}},\quad(\dphi\lam)a=\bw_{b\in B}\lam b\rat\phi(a,b),\label{dphi-tensor-cotensor-def}
    \end{align}
    which satisfies
    \begin{equation} \label{phi-tensor-cotensor-Isbell}
    X^{A^{\op}}(\mu,\dphi\lam)=(\sV^B)^{\op}(\uphi\mu,\lam)=\bw_{b\in B}\lam b\ra X^{A^{\op}}(\mu,\phi(-,b))
    \end{equation}
    for all $\mu\in X^{A^{\op}}$, $\lam\in\sV^B$. In particular:
    \begin{enumerate}[label=(1\alph*)]
    \item If $B={\bf 1}$, the singleton $\sV$-category, then a $\sV$-bifunctor $\tau\colon A^{\op}\ts{\bf 1}\to X$ may be identified with $\tau\in X^{A^{\op}}$, and \eqref{phi-tensor-cotensor-Isbell} becomes
        $$X^{A^{\op}}(\mu,\dtau v)=\sV^{\op}(\utau\mu,v)=v\ra X^{A^{\op}}(\mu,\tau)$$
        for all $\mu\in X^{A^{\op}}$, $v\in\sV$; that is,
        $$\dtau v=v\rat\tau$$
        is the cotensor of $v$ and $\tau$ in $X^{A^{\op}}$ (cf. \eqref{tensor-def}), given by the pointwise cotensors inherited from $X$, and
        $$\sM\tau=\{v\rat\tau\mid v\in\sV\}.$$
    \item If
        \begin{equation} \label{eval-def}
        \phi\colon A^{\op}\ts X^{A^{\op}}\to X,\quad\phi(a,\mu)=\mu a
        \end{equation}
        is the \emph{evaluation} $\sV$-bifunctor, then it follows from \eqref{uphi-tensor-cotensor-def} and \eqref{dphi-tensor-cotensor-def} that
        $$\uphi\mu=X^{A^{\op}}(\mu,-)\quad\text{and}\quad\dphi\Lam=\bw_{\tau\in X^{A^{\op}}}\Lam\tau\rat\tau$$
        for all $\mu\in X^{A^{\op}}$, $\Lam\in\sV^{X^{A^{\op}}}$; that is,
        $$\uphi=\syd_{X^{A^{\op}}}\quad\text{and}\quad\dphi={\inf}_{X^{A^{\op}}},$$
        where the latter identity is an immediate consequence of the completeness of $X^{A^{\op}}$ and the formula for infima given in \cite[Theorem 2.8]{Shen2013a}. Since $\dphi\uphi={\inf}_{X^{A^{\op}}}\syd_{X^{A^{\op}}}\cong 1_{X^{A^{\op}}}$, we have now
        $$\Mphi=X^{A^{\op}}.$$
    \end{enumerate}
\item The Kan adjunction induced by a $\sV$-bifunctor $\psi\colon A^{\op}\ts B\to X$ is given by
    \begin{align}
    &\psi^*\colon \sV^{B^{\op}}\to X^{A^{\op}},\quad(\psi^*\lam)a=\bv_{b\in B}\lam b\star\psi(a,b),\label{psi-ustar-tensor-cotensor-def}\\
    &\psi_*\colon X^{A^{\op}}\to\sV^{B^{\op}},\quad(\psi_*\mu)b=\bw_{a\in A}X(\psi(a,b),\mu a)=X^{A^{\op}}(\psi(-,b),\mu),\label{psi-lstar-tensor-cotensor-def}
    \end{align}
    which satisfies
    \begin{equation} \label{psi-tensor-cotensor-Kan}
    X^{A^{\op}}(\psi^*\lam,\mu)=\sV^{B^{\op}}(\lam,\psi_*\mu)=\bw_{b\in B}\lam b\ra X^{A^{\op}}(\psi(-,b),\mu)
    \end{equation}
    for all $\mu\in X^{A^{\op}}$, $\lam\in\sV^{B^{\op}}$. In particular:
    \begin{enumerate}[label=(2\alph*)]
    \item If $B={\bf 1}$, then for each $\tau\in X^{A^{\op}}$, \eqref{psi-tensor-cotensor-Kan} becomes
        $$X^{A^{\op}}(\tau^*v,\mu)=\sV(v,\tau_*\mu)=v\ra X^{A^{\op}}(\tau,\mu)$$
        for all $\mu\in X^{A^{\op}}$, $v\in\sV$; that is,
        $$\tau^*v=v\star\tau$$
        is the tensor of $v$ and $\tau$ in $X^{A^{\op}}$ (cf. \eqref{tensor-def}), given by the pointwise tensors inherited from $X$, and
        $$\sK\tau=\{X^{A^{\op}}(\tau,\mu)\mid\mu\in X^{A^{\op}}\}.$$
    \item If $\psi\colon A^{\op}\ts X^{A^{\op}}\to X$ is the evaluation $\sV$-bifunctor given by \eqref{eval-def}, then it follows from \eqref{psi-ustar-tensor-cotensor-def} and \eqref{psi-lstar-tensor-cotensor-def} that
        $$\psi^*\Lam=\bv_{\tau\in X^{A^{\op}}}\Lam\tau\star\tau\quad\text{and}\quad\psi_*\mu=X^{A^{\op}}(-,\mu)$$
        for all $\mu\in X^{A^{\op}}$, $\Lam\in\sV^{(X^{A^{\op}})^{\op}}$; that is,
        $$\psi^*={\sup}_{X^{A^{\op}}}\quad\text{and}\quad\psi_*=\sy_{X^{A^{\op}}}$$
        by \cite[Theorem 2.8]{Shen2013a}, and consequently
        $$\Kpsi=\Fix(\sy_{X^{A^{\op}}}{\sup}_{X^{A^{\op}}})=\{\sy_{X^{A^{\op}}}\mu\mid\mu\in X^{A^{\op}}\}\simeq X^{A^{\op}}.$$
    \end{enumerate}
\item The Kan adjunction induced by a $\sV$-bifunctor $\zeta\colon A^{\op}\ts B\to\sV$, i.e., a $\sV$-distributor $\zeta\colon A\oto B$, is given by
    \begin{align}
    &\zeta_{\dag}\colon (X^B)^{\op}\to (X^A)^{\op},\quad(\zeta_{\dag}\lam)a=\bw_{b\in B}\zeta(a,b)\rat\lam b, \label{zeta-ldag-tensor-cotensor-def}\\
    &\zeta^{\dag}\colon(X^A)^{\op}\to(X^B)^{\op},\quad(\zeta^{\dag}\mu)b=\bv_{a\in A}\zeta(a,b)\star\mu a. \label{zeta-udag-tensor-cotensor-def}
    \end{align}
    Note that $\zeta_{\dag}\dv\zeta^{\dag}$ apparently generalizes the adjunction given by \cite[Proposition 6.2.1]{Shen2014} (in the case of $\CQ=\sV$) by replacing the multiplication $\otimes$ and the implication $\ra$ of $\sV$ with the tensor and cotensor of $X$, respectively. In particular, if $A=(A,\al)$ is a $\sV$-category, then $\al\colon A\oto A$ is a $\sV$-distributor, and it is straightforward to check that $\al_{\dag}\cong\al^{\dag}\cong 1_{(X^A)^{\op}}$, which necessarily forces that $\sK^{\dag}\al=(X^A)^{\op}$.
\end{enumerate}


\subsection{(Generalized) metric spaces}

Recall that a (generalized) metric space is precisely a $[0,\infty]$-category (cf. Example \ref{V-cat-exmp}\ref{V-cat-exmp:metric}). The implication of the quantale
$$[0,\infty]=([0,\infty],\geq,+,0)$$
is given by
$$x\ra y=\max\{0,y-x\}$$
for all $x,y\in[0,\infty]$, where ``$\max$'' refers to the standard order of extended real numbers. In order to eliminate ambiguity, we denote joins and meets with respect to the standard order on $[0,\infty]$ by $\bigsqcup$ and $\bigsqcap$, respectively.

Given a topological space $X$, let $\LSC$ denote the set of lower semi-continuous maps from $X$ to $[0,\infty]$; that is, maps $f\colon X\to[0,\infty]$ such that $\{x\in X\mid fx\leq a\}$ is closed for all $a\in[0,\infty]$. $\LSC$ becomes a (generalized) metric space, with
$$\LSC(f,g)=\bigsqcup_{x\in X}\max\{0,gx-fx\}$$
for all $f,g\in\LSC$. It is clear that $\LSC$ is closed with respect to pointwise joins of maps, making itself a complete lattice, and moreover:
\begin{enumerate}[label=(\arabic*)]
\item $\LSC$ is tensored, with
    $$a\star f\colon X\to[0,\infty],\quad (a\star f)x=a+fx$$
    for all $a\in[0,\infty]$, $f\in\LSC$, $x\in X$;
\item $\LSC$ is cotensored, with
    $$a\rat f\colon X\to[0,\infty],\quad (a\rat f)x=\max\{0,fx-a\}$$
    for all $a\in[0,\infty]$, $f\in\LSC$, $x\in X$.
\end{enumerate}
It follows from Proposition \ref{V-complete-tensor} that $\LSC$ is a complete $[0,\infty]$-category. Hence, by setting $\sV=[0,\infty]$ and $X=\LSC$ in \eqref{tensor-cotensor-two-variable} we obtain a two-variable adjunction in metric spaces.

Let $Y$ be a classical metric space (whose metric is symmetric, finite and separated) equipped with the usual metric topology. Then the set $\sC(X,Y)$ of continuous maps from $X$ to $Y$, and the set $[0,\infty]^Y$ of non-expansive maps from $Y$ to $[0,\infty]$, are both equipped with a (generalized) metric structure given by Equation \eqref{Y-X-f-g} in Example \ref{V-cat-exmp}\ref{V-cat-exmp:functor}. Since the composition $gf$ of $f\in\sC(X,Y)$ and $g\in[0,\infty]^Y$ lies in $\LSC$, it is straightforward to check that
$$\phi\colon[0,\infty]^Y\times\sC(X,Y)\to\LSC,\quad (g,f)\mapsto gf$$
defines a $[0,\infty]$-bifunctor, where the metric on $[0,\infty]^Y\times\sC(X,Y)$ is given by Example \ref{V-cat-exmp}\ref{V-cat-exmp:ts}.
\begin{enumerate}[label=(\arabic*)]
\item The induced Isbell adjunction
    $$\uphi\dv\dphi\colon([0,\infty]^{\sC(X,Y)})^{\op}\to \LSC^{[0,\infty]^Y}$$
    is given by
    \begin{align*}
    &(\uphi\mu)f=\bigsqcup_{g\in[0,\infty]^Y}\LSC(\mu g,gf)=\bigsqcup_{g\in[0,\infty]^Y}\bigsqcup_{x\in X}\max\{0,gfx-(\mu g)x\},\\
    &((\dphi\lam)g)x=\bigsqcup_{f\in\sC(X,Y)}(\lam f\rat gf)x=\bigsqcup_{f\in\sC(X,Y)}\max\{0,gfx-\lam f\}
    \end{align*}
    for all $\mu\in\LSC^{[0,\infty]^Y}$, $\lam\in[0,\infty]^{\sC(X,Y)}$, and $\Mphi$ consists of $\mu\in\LSC^{[0,\infty]^Y}$ satisfying
    $$(\mu g)x=\bigsqcup_{f\in\sC(X,Y)}\bigsqcap_{g'\in[0,\infty]^Y}\bigsqcap_{x'\in X}\max\{0,gfx-\max\{0,g'fx'-(\mu g')x'\}\}$$
    for all $g\in[0,\infty]^Y$, $x\in X$. 
\item The induced Kan adjunction
    $$\phi^*\dv\phi_*\colon \LSC^{[0,\infty]^Y}\to[0,\infty]^{\sC(X,Y)^{\op}}$$
    is given by
    \begin{align*}
    &((\phi^*\lam)g)x=\bigsqcap_{f\in\sC(X,Y)}(\lam f\star gf)x=\bigsqcap_{f\in\sC(X,Y)}\lam f+gfx,\\
    &(\phi_*\mu)f=\bigsqcup_{g\in[0,\infty]^Y}\LSC(gf,\mu g)=\bigsqcup_{g\in[0,\infty]^Y}\bigsqcup_{x\in X}\max\{0,(\mu g)x-gfx\}
    \end{align*}
    for all $\lam\in[0,\infty]^{\sC(X,Y)^{\op}}$, $\mu\in\LSC^{[0,\infty]^Y}$, and $\Kphi$ consists of $\lam\in[0,\infty]^{\sC(X,Y)^{\op}}$ satisfying
    $$\lam f=\bigsqcup_{g\in[0,\infty]^Y}\bigsqcup_{x\in X}\bigsqcap_{f'\in\sC(X,Y)}\max\{0,\lam f'+gf'x-gfx\}$$
    for all $f\in\sC(X,Y)$.
\item In particular, for each $x\in X$,
    $$\phi_x\colon[0,\infty]^Y\times\sC(X,Y)\to[0,\infty],\quad (g,f)\mapsto gfx$$
    defines a $[0,\infty]$-distributor. The induced Kan adjunction
    $$(\phi_x)_{\dag}\dv\phi_x^{\dag}\colon(\LSC^{([0,\infty]^Y)^{\op}})^{\op}\to(\LSC^{\sC(X,Y)})^{\op}$$
    is given by
    \begin{align*}
    &(((\phi_x)_{\dag}\lam)g)x'=\bigsqcup_{f\in\sC(X,Y)}(gfx\rat\lam f)x'=\bigsqcup_{f\in\sC(X,Y)}\max\{0,\lam fx'-gfx\}, \\
    &((\phi_x^{\dag}\mu)f)x'=\bigsqcap_{g\in[0,\infty]^Y}(gfx\star\mu g)x'=\bigsqcap_{g\in[0,\infty]^Y}gfx+(\mu g)x'.
    \end{align*}
    for all $\lam\in\LSC^{\sC(X,Y)}$, $\mu\in\LSC^{([0,\infty]^Y)^{\op}}$, and $\Kdphi$ consists of $\lam\in\LSC^{\sC(X,Y)}$ satisfying
    $$(\lam f)x'=\bigsqcap_{g\in[0,\infty]^Y}\bigsqcup_{f'\in\sC(X,Y)}(gfx+\max\{0,\lam f'x'-gf'x\})$$
    for all $f\in\sC(X,Y)$, $x'\in X$.
\end{enumerate}

\begin{acknowledgements}
This work was supported by National Natural Science Foundation of China (No. 12071319). The authors would like to thank Professor Hongliang Lai and Professor Dexue Zhang for helpful discussions. The authors would also like to thank the anonymous referee for valuable comments and suggestions which significantly improve the presentation of the paper.
\end{acknowledgements}

%
\section*{Conflict of interest}
The authors declare that they have no conflict of interest.

\section*{Data availability statement}
Data sharing not applicable to this article as no datasets were generated or analysed during the current study.



%
%

\end{document}